\documentclass[12pt]{article}
\usepackage{amsmath,amsthm,amsfonts,amssymb,mathrsfs,latexsym,graphicx}
\newtheorem{thm}{Theorem}[section]
\newtheorem{lem}{Lemma}[section]
\newtheorem{cor}{Corollary}[section]
\newtheorem{conj}{Conjecture}[section]
\newtheorem{rem}{Remark}[section]
\newtheorem{prop}{Proposition}[section]

\newcommand{\lrf}[1]{\left\lfloor #1\right\rfloor}
\newcommand{\lrc}[1]{\left\lceil #1\right\rceil}

\begin{document}
\title{Proof of a conjecture on unimodality}
\author{
Yi\ Wang$^{\rm a}$ \footnote{\footnotesize{\it E-mail addresses:}
wangyi@dlut.edu.cn
(Y. Wang), mayeh@math.sinica.edu.tw (Y.-N. Yeh)},\quad Yeong-Nan\ Yeh$^{\rm b}$\\
{\footnotesize\it $^{\rm a}$ Department of Applied Mathematics,
Dalian University
of Technology, Dalian 116024, P. R. China}\\[5pt]
{\footnotesize\it $^{\rm b}$ Institute of Mathematics, Academia
Sinica, Taipei 11529, Taiwan}
}

\date{\small Received 30 September 2003;\quad accepted 26 April 2004}

\maketitle

\begin{abstract}Let $P(x)$ be a polynomial of degree $m$, with
nonnegative and non-decreasing coefficients. We settle the
conjecture that for any positive real number $d$, the coefficients
of $P(x+d)$ form a unimodal sequence, of which the special case $d$
being a positive integer has already been asserted in a previous
work. Further, we explore the location of modes of $P(x+d)$ and
present some sufficient conditions on $m$ and $d$ for which $P(x+d)$
has the unique mode $\left\lceil{m-d\over d+1}\right\rceil$.
\\[7pt]
{\it MSC:}\quad 05A20\\[7pt]
{\it Keywords:}\quad Log-concavity; Unimodality; Modes
\end{abstract}

\section{Introduction}
\hspace*{\parindent}
Let $a_0,\ldots,a_m$ be a sequence of nonnegative real numbers.
We say that the sequence is {\em unimodal}
if there exists an index $0\le t\le m$
such that $a_0\le \cdots\le a_{t-1}\le a_t\ge a_{t+1}\ge\cdots\ge a_m$.
Such an index $t$ is called a mode of the sequence.
A property closely related to unimodality is log-concavity.
We say that the sequence is {\em log-concave}
if $a_{i-1}a_{i+1}\le a_i^2$ for all $1\le i\le m-1$.
The sequence is said to {\em have no internal zeros}
if there are not three indices $i<j<k$ such that
$a_i,a_k\not=0$ and $a_j=0$.
It is well known that a log-concave sequence
with no internal zeros is unimodal
(see \cite[Proposition 2.5.1]{Bre89} for instance).
Unimodal and log-concave sequences occur naturally
in many branches of mathematics.
See the survey articles \cite{Sta89} and \cite{Bre94}
for various techniques, problems, and results
about unimodality and log-concavity.

Let $P(x)=\sum\limits_{i=0}^ma_ix^i$ be a polynomial with nonnegative coefficients.
We say that $P(x)$ is unimodal
(respectively, log-concave, non-decreasing, etc.)
if the sequence of coefficients $a_0,a_1,\ldots,a_m$ of $P(x)$
enjoys the corresponding property.
A mode of $a_0,\ldots,a_m$ is also called a mode of $P(x)$.

It is well known that if $P(x)$ is log-concave with no internal zeros,
then $P(x+1)$ is log-concave
(see \cite[Corollary 8.4]{Bre94} or \cite[Theorem 2]{Hog74}).
Actually, it may also be shown that
$P(x+d)$ is log-concave for any positive number $d$ by using \cite[Theorem 2.5.3]{Bre89}.
In the present paper,
we consider the analogue problem concerning unimodality.
Let $P(x)$ be nonnegative and non-decreasing.
It is shown that $P(x+1)$ is unimodal in \cite{BM99}
and more generally,
that $P(x+n)$ is unimodal when $n$ is a positive integer in \cite{AABKMR01}.
Further, the following is conjectured.
\begin{conj}[\cite{AABKMR01}]\label{conj1}
Let $P(x)$ be a polynomial of degree $m$ and with nonnegative
coefficients. Suppose that $P(x)$ is non-decreasing and that $d$
is a positive real number. Then $P(x+d)$ is unimodal.
\end{conj}

In this paper we settle the above conjecture.
Moreover,
we will explore the number and location of modes of the polynomial $P(x+d)$.
Let $M_*(P,d)$ and $M^*(P,d)$ be
the smallest and the greatest mode of $P(x+d)$ respectively.
Denote ${\overline m}(d)=\lrc{{m-d \over d+1}}$
and ${\underline m}(d)=\lrf{{m\over d+1}}$
where $\lrc x$ and $\lrf x$ denote the least integer $\ge x$
and the greatest integer $\le x$ respectively.
It is not difficult to see that
${\overline m}(d)$ and ${\underline m}(d)$ coincide when $d$ is a positive integer.
In \cite{AABKMR01},
it is shown that ${\underline m}(d)$ is a mode of $P(x+d)$
when $d$ is a positive integer.
The statement is not true when $d$ is only a positive number.
Generally speaking,
the number and location of modes of $P(x+d)$
are related not only to $m$ and $d$,
but also to coefficients of the polynomial $P(x)$.
The matter is somewhat different when $d\ge 1$.
In this case, we can show that
$P(x+d)$ has at most two modes
${\overline m}(d)$ and ${\overline m}(d)+1$
if $P(x)=ax^m$,
or ${\overline m}(d)-1$ and ${\overline m}(d)$ otherwise.
We will also present certain sufficient conditions on $m$ and $d$
that $P(x+d)$ has the unique mode ${\overline m}(d)$,
including the case when $d$ is a positive integer larger than $1$.

Throughout this paper,
let $m$ be a positive integer and $d$ a positive real number.
We denote by ${\bf P}^m_\uparrow$ the set of monic polynomials of degree $m$,
with nonnegative and non-decreasing coefficients.
When there is no danger of confusion,
we simply write ${\overline m}$ to mean ${\overline m}(d)$.
By definition, it follows immediately that
\begin{eqnarray}\label{basic}
m-d\le (d+1){\overline m}<m+1,
\end{eqnarray}
which will be used repeatedly in the sequel.

\section{Proof of Conjecture \ref{conj1}}
\hspace*{\parindent}
To prove Conjecture \ref{conj1},
we need the following two lemmas.
\begin{lem}\label{lem2.1}
Suppose that the polynomial $f(x)$ is unimodal with the smallest mode $t$
and that $d>0$.
Then $(x+d)f(x)$ is unimodal with the smallest mode $t$ or $t+1$.
\end{lem}
\proof Let $f(x)=\sum\limits_{i=0}^nc_ix^i$ where $c_0\le\cdots\le
c_{t-1}<c_t\ge c_{t+1}\ge \cdots\ge c_n$. Then
\begin{eqnarray*} (x+d)f(x)&=&c_0d+(c_0+c_1d)x+\cdots+(c_{t-2}+c_{t-1}d)x^{t-1}+(c_{t-1}+c_td)x^t\\
&&+(c_{t}+c_{t+1}d)x^{t+1}+\cdots+(c_{n-1}+c_nd)x^{n-1}+c_nx^n.
\end{eqnarray*}
Clearly, $c_0\le c_0+c_1d\le\cdots\le c_{t-2}+c_{t-1}d<c_{t-1}+c_td$
and $c_t+c_{t+1}d\ge \cdots\ge c_{n-1}+c_nd\ge c_n$. So the
statement follows. \qed
\begin{lem}\label{lem2.2}
Let $P(x)=\sum\limits_{i=0}^ma_ix^i$ be a polynomial of degree $m$,
with nonnegative coefficients and $d>0$.
Suppose that $P(x+d)=\sum\limits_{j=0}^mb_jx^j$.
Then $b_{\overline m}\ge b_{{\overline m}+1}\ge\cdots\ge b_m$.
Furthermore, if $d\ge (m-1)/2$, then $P(x+d)$ is unimodal and has the mode $0$ or $1$.
In particular, if $d\ge m$ then $P(x+d)$ is non-increasing.
\end{lem}
\proof We have
$b_j={P^{(j)}(d)/j!}=\sum\limits_{i=j}^ma_id^{i-j}{i\choose j}$,
which yields that
\begin{eqnarray}\label{eq1}
(j+1)d^{j+1}(b_{j+1}-b_{j})
=\sum\limits_{i=j}^ma_id^{i}{i\choose j}[(i+1)-(d+1)(j+1)].
\end{eqnarray}
Now let $j\ge {\overline m}$. Then $(d+1)(j+1)\ge (d+1)({\overline
m}+1)\ge m+1$ by (\ref{basic}). Every term in the sum (\ref{eq1}) is
therefore non-positive, and thus $b_{j+1}\le b_{j}$. Finally, note
that $(m-1)/2\le d<m$ implies ${\overline m}\le 1$, and that $d\ge
m$ implies ${\overline m}=0$. So the statement follows. \qed

\proof[Proof of Conjecture \ref{conj1}] Let
$P(x)=\sum\limits_{i=0}^ma_ix^i$ and
$P(x+d)=\sum\limits_{j=0}^mb_jx^j$. We need to show that
$b_0,\ldots,b_m$ is unimodal. We do this by induction on $m$. If
$m=1$, the result is obvious, so we proceed to the inductive step.
By Lemma \ref{lem2.2}, it suffices to consider the case $m>2d+1$.

Let $P(x)=a_0+xf(x)$ where $f(x)=\sum\limits_{i=0}^{m-1}a_{i+1}x^i$.
Then
$$P(x+d)=a_0+(x+d)f(x+d).$$
By the induction hypothesis,
$f(x+d)$ is unimodal,
so is $(x+d)f(x+d)$ by Lemma \ref{lem2.1}.
Thus $b_1,b_2,\ldots,b_m$ is
unimodal.

Let $r=\lrf{d}$. Then $r<d+1<m$.
By (\ref{eq1}) we have
\begin{eqnarray*} b_1-b_0
&=&\sum\limits_{i=0}^ma_id^{i-1}(i-d)\\
&=&\sum\limits_{i=r+1}^ma_id^{i-1}(i-d)-\sum\limits_{i=0}^{r}a_id^{i-1}(d-i)\\
&\ge&a_r\sum\limits_{i=r+1}^md^{i-1}(i-d)-a_r\sum\limits_{i=0}^{r}d^{i-1}(d-i)\\
&=&a_r[d+2d^2+\cdots+(m-1)d^{m-1}-d^m]\\
&\ge&a_r[(m-1)d^{m-1}-d^m]\\
&=&a_r(m-d-1)d^{m-1}\\
&\ge& 0.
\end{eqnarray*}
Thus $b_0,b_1,\ldots,b_m$ is still unimodal. This completes the
proof. \qed
\begin{cor}\label{cor2.1}
Let $P(x)\in{\bf P}^m_\uparrow$ and $d>0$.
Suppose that $P(x)\ne x^m$.
Then $$M^*(P,d)\le {\overline m}.$$
\end{cor}
\proof Let $P(x)=\sum\limits_{i=0}^ma_ix^i$ and
$P(x+d)=\sum\limits_{j=0}^mb_jx^j$. We have by (\ref{eq1})
$$
({\overline m}+1)d^{{\overline m}+1}
(b_{{\overline m}+1}-b_{\overline m})
=\sum\limits_{i={\overline m}}^m
a_id^{i}{i\choose {\overline m}}[(i+1)-(d+1)({\overline m}+1)].
$$
By (\ref{basic}), $(i+1)-(d+1)({\overline m}+1) \le
(m+1)-(d+1)({\overline m}+1)\le 0$ for each $i\le m$. In particular,
$m-(d+1)({\overline m}+1)\le -1<0$. On the other hand, $a_{m-1}\ne
0$ since $P(x)\ne x^m$. Hence $b_{{\overline m}+1}<b_{\overline m}$.
This implies that the unimodal sequence $\{b_j\}$ has no mode larger
than ${\overline m}$, and the proof is therefore complete. \qed

\section{Modes of $(x+d)^m$ and $\sum\limits_{i=0}^m(x+d)^i$}
\hspace*{\parindent}
This section is devoted to studying
modes of $P(x+d)$ for two basic polynomials
$P(x)=x^m$ and $P(x)=\sum\limits_{i=0}^mx^i$ respectively,
which will play a key role in investigating modes of $P(x+d)$
for generic polynomials $P(x)\in{\bf P}^m_\uparrow$.
\begin{prop}\label{x^m}
Let $d>0$.
If ${m+1\over d+1}\in{\bf Z}^+$,
then $(x+d)^m$ has two modes
${\overline m}$ and ${\overline m}+1$;
otherwise $(x+d)^m$ has the unique mode ${\overline m}$.
\end{prop}
\proof Let $(x+d)^m=\sum\limits_{i=0}^mc_ix^i$ where $c_i={m\choose
i}d^{m-i}$. Denote $f(x)={m-x+1\over dx}$. Then ${c_i\over
c_{i-1}}=f(i)$. Clearly, $f(x)$ is strictly decreasing and
$f({m+1\over d+1})=1$. Now $i\le {\overline m}$ implies $i<{m+1\over
d+1}$, and $i\ge {\overline m}+1$ implies $i\ge {m+1\over d+1}$. So
the statement follows. \qed

Let $Q_m(x)=\sum\limits_{i=0}^mx^i$
and $Q_m(x+d)=\sum\limits_{j=0}^md_jx^j$ where
\begin{eqnarray}\label{qmd}
d_j=\sum\limits_{i=j}^md^{i-j}{i\choose j},\qquad j=0,1,\ldots,m.
\end{eqnarray}
Then the sequence $\{d_j\}$ is log-concave with no internal
zeros(see Brenti\cite[Theorem 2.5.3]{Bre89} for instance).
Actually, we have the following stronger result.
\begin{prop}
The sequence $\{d_j\}$ is strictly log-concave,
i.e., $d_{j-1}d_{j+1}<d_j^2$ for all $0<j<m$,
and is therefore unimodal with at most two modes.
\end{prop}
\proof Note that
\begin{eqnarray*}
d_{j-1}&=&\sum\limits_{i=j-1}^md^{i-j+1}{i\choose j-1}\\
&=&\sum\limits_{i=j-1}^md^{i-j+1}\left[{i+1\choose j}-{i\choose j}\right]\\
&=&(1-d)d_j+d^{m-j+1}{m+1\choose j}.
\end{eqnarray*}
Thus we have
\begin{eqnarray*}
d_j^2-d_{j-1}d_{j+1}
&=&d_j^2-\left[(1-d)d_j+d^{m-j+1}{m+1\choose j}\right]d_{j+1}\\
&=&\left[d_j-(1-d)d_{j+1}\right]d_j-d^{m-j+1}{m+1\choose j}d_{j+1}\\
&=&d^{m-j}{m+1\choose j+1}d_j-d^{m-j+1}{m+1\choose j}d_{j+1}\\
&=&\sum_{i=j}^m
\left[{m+1\choose j+1}{i\choose j}-{m+1\choose j}{i\choose j+1}\right]d^{m+i-2j}\\
&=&\sum_{i=j}^m
{m-i+1\over j+1}{m+1\choose j}{i\choose j}d^{m+i-2j}\\
&>&0,
\end{eqnarray*}
the desired inequality. \qed

In what follows we explore the location of modes of the sequence $\{d_j\}$.
We first consider the case $d\ge 1$.
The matter is rather simple when $d=1$.
\begin{prop}\label{Q_m,1}
If $m$ is even then $Q_m(x+1)$ has two modes ${m\over2}-1$ and ${m\over2}$;
otherwise $Q_m(x+1)$ has the unique mode ${m-1\over2}$.
\end{prop}
\proof Since $Q_m(x)={1\over x-1}(x^{m+1}-1)$, we have
$$Q_m(x+1)={1\over x}\left[(x+1)^{m+1}-1\right].$$
By Proposition \ref{x^m}, $(x+1)^{m+1}$ has two modes $\overline
{m+1}={m\over2}$ and $\overline {m+1}+1={m\over2}+1$ for $m$ even,
or only one mode $\overline {m+1}={m+1\over2}$ otherwise, so does
$(x+1)^{m+1}-1$. Thus the statement follows. \qed
\begin{prop}\label{Q_m,d>1}
Let $d\ge 1$.
Then $Q_m(x+d)$ has at most two modes
$\overline {m}-1$ and $\overline {m}$.
In particular, if $\overline {m+1}={\overline m}+1$,
then $Q_m(x+d)$ has the unique mode ${\overline m}$.
\end{prop}
\proof By Lemma \ref{lem2.2}, it suffices to consider the case $1\le
d<m$. We have
$$(x+d-1)Q_m(x+d)=(x+d)^{m+1}-1.$$
By Proposition \ref{x^m},
$(x+d)^{m+1}$ has the smallest mode $\overline {m+1}$,
so does $(x+d)^{m+1}-1$.
Thus $M_*(Q_m,d)\ge \overline {m+1}-1$ by Lemma \ref{lem2.1}.
On the other hand,
we have $M^*(Q_m,d)\le {\overline m}$ by Corollary \ref{cor2.1}.
Note that $\overline {m+1}={\overline m}$ or ${\overline m}+1$ since
$${m-d\over d+1}<{m+1-d\over d+1}<{m-d\over d+1}+1.$$
Hence $Q_m(x+d)$ has at most two modes ${\overline m}-1$ and
${\overline m}$, and in particular, only one mode ${\overline m}$
provided $\overline {m+1}={\overline m}+1$. This completes the
proof. \qed
\begin{cor}\label{integ1}
If $d\ge 1$ and ${m+1\over d+1}\in{\bf Z}^+$, then $Q_m(x+d)$ has
the unique mode ${\overline m}$.
\end{cor}
\proof If ${m+1\over d+1}\in{\bf Z}^+$, then ${m-d\over d+1}\in{\bf
Z}^+$, and so ${\overline m}={m-d\over d+1}$. On the other hand,
$$\overline {m+1}=\lrc{{m+1-d\over d+1}}=
\lrc{{m+1\over d+1}-{d\over d+1}} ={m+1\over d+1}.$$ Thus $\overline
{m+1}={\overline m}+1$. So the statement follows from Proposition
\ref{Q_m,d>1}. \qed
\begin{prop}\label{d>1,1}
If $d>1$ and $d{\overline m}\in{\bf Z}^+$,
then $Q_m(x+d)$ has the unique mode ${\overline m}$.
\end{prop}
\proof By Proposition \ref{Q_m,d>1}, it suffices to prove
$d_{{\overline m}}>d_{{\overline m}-1}$. By (\ref{eq1}), we have
$${\overline m}d^{{\overline m}}(d_{{\overline m}}-d_{{\overline m}-1})
=\sum_{i={\overline m}-1}^md^i{i\choose {\overline m}-1}[(i+1)-(d+1){\overline m}].$$
The sum contains terms of both signs.
Let $r=\lrc{(d+1){\overline m}}-1$. Denote
$$S_1=\sum_{i=r}^md^i{i\choose {\overline m}-1}[(i+1)-(d+1){\overline m}]$$
and
$$
S_2=\sum_{i={\overline m}-1}^{r-1}d^i
{i\choose {\overline m}-1}[(d+1){\overline m}-(i+1)].
$$
Then
${\overline m}d^{{\overline m}} (d_{{\overline m}}-d_{{\overline m}-1})=S_1-S_2$.
Thus we need to prove $S_1>S_2$.

Since $(d+1){\overline m}<m+1$ by (\ref{basic})
and the left is an integer by the assumption,
we have $r\le m-1$. So
$$
S_1\ge d^{r+1}{r+1\choose {\overline m}-1}[(r+2)-(d+1){\overline m}]
= d^{r+1}{r+1\choose {\overline m}-1}.
$$
On the other hand,
\begin{eqnarray*}
S_2
&\le&
\sum\limits_{i={\overline m}-1}^{r-1}d^{r-1}
{i\choose {\overline m}-1}[(r+1)-(i+1)]\\
&\le&
d^{r-1}\left[(r+1)\sum\limits_{i={\overline m}-1}^{r-1}
{i\choose {\overline m}-1}
-{\overline m}\sum\limits_{i={\overline m}-1}^{r-1}
{i+1\choose {\overline m}-1}\right]\\
&=&
d^{r-1}\left[(r+1){r\choose {\overline m}}
-{\overline m}{r+1\choose {\overline m}+1}\right]\\
&=&
d^{r-1}{r+1\choose {\overline m}+1}.
\end{eqnarray*}
Thus we have
\begin{eqnarray*}
{S_1\over S_2}
\ge {d^{r+1}{r+1\choose {\overline m}-1}\over d^{r-1}{r+1\choose {\overline m}+1}}
= {d^2{\overline m}({\overline m}+1)\over (r-{\overline m}+1)(r-{\overline m}+2)}
= {d({\overline m}+1)\over d{\overline m}+1}>1,
\end{eqnarray*}
the desired inequality. \qed
\begin{cor}
If $d>1$ and $d\in {\bf Z}^+$,
then $Q_m(x+d)$ has the unique mode ${\overline m}$.
\end{cor}
\begin{cor}
If $d>1$ and ${m\over d+1}\in{\bf Z}^+$,
then $Q_m(x+d)$ has the unique mode ${\overline m}$.
\end{cor}
\proof If ${m\over d+1}\in{\bf Z}^+$, then
$${\overline m}=\lrc{m-d\over d+1}
=\lrc{{m\over d+1}-{d\over d+1}}={m\over d+1}.$$ Thus $d{\overline
m}=m-{\overline m}\in{\bf Z}^+$, and the statement follows from
Proposition \ref{d>1,1}. \qed

We next consider the case $0<d<1$,
which is more complicated.
For example,
modes of $Q_m(x+d)$ may be neither ${\overline m}-1$ nor ${\overline m}$
(see Remark \ref{rem3.1}).
The following is some rough estimate
for location of modes of $Q_m(x+d)$.
\begin{prop}\label{Q_m,0<d<1}
Let $0<d<1$. Then

{\em (i)}\quad
$\lrf{{m\over2}}\le M_*(Q_m,d)\le M^*(Q_m,d)
\le \min\{m-1,{\overline m}\}$.

{\em (ii)}\quad
If $0<d<1/{m\choose 2}$,
then $Q_m(x+d)$ has the unique mode $m-1$.
The converse is also true.

{\em (iii)}\quad
If $0<1-d\le 1/m$,
then $Q_m(x+d)$ has at most two modes ${\overline m}-1$ and ${\overline m}$.
In particular, if ${m+1\over d+1}\in{\bf Z}^+$,
then $Q_m(x+d)$ has the unique mode ${\overline m}$.

{\em (iv)}\quad
There exists a positive number $\varepsilon$ such that for $0<1-d<\varepsilon$,
$Q_m(x+d)$ has the unique mode $\lrf{m\over2}$.
\end{prop}
\proof By the definition, $M_*(Q_m,d)$ is the greatest integer $j$
no larger than ${\overline m}$ such that $d_j>d_{j-1}$. Note that
\begin{eqnarray}\label{dine}
d_{j}-d_{j-1}
&=&\sum\limits_{i=j}^m{i\choose j}d^{i-j}-
\sum\limits_{i=j-1}^m{i\choose j-1}d^{i-j+1}\nonumber\\
&=&\sum\limits_{i=j-1}^{m-1}{i+1\choose j}d^{i-j+1}-
\sum\limits_{i=j-1}^m{i\choose j-1}d^{i-j+1}\nonumber\\
&=&\sum\limits_{i=j}^{m-1}{i\choose j}d^{i-j+1}-{m\choose j-1}d^{m-j+1}.
\end{eqnarray}
Hence
\begin{eqnarray*}
M_*(Q_m,d)
=\max\left\{1\le j\le {\overline m}:
\sum\limits_{i=j}^{m-1}{i\choose j}d^{i-j+1}
-{m\choose j-1}d^{m-j+1}>0\right\}.
\end{eqnarray*}
When $0<d<1$, we have
\begin{eqnarray*}
\sum\limits_{i=j}^{m-1}{i\choose j}d^{i-j+1}
\ge d^{m-j}\sum\limits_{i=j}^{m-1}{i\choose j}
={m\choose j+1}d^{m-j}.
\end{eqnarray*}
It is not difficult to see that
$${m\choose j+1}d^{m-j-1}-{m\choose j-1}d^{m-j}>0$$
is equivalent to
$$(m-j)(m-j+1)-dj(j+1)>0.$$
Now let $h(x)=(m-x)(m-x+1)-dx(x+1)$.
Then $h(x)$ is a decreasing function in the interval $0\le x\le m$ since $h'(x)<0$.
Thus $h(x_0)>0$ for some $x_0\in (0,m)$
implies that $M_*(Q_m,d)\ge\lrf{x_0}$.

(i)\quad
Since
$h\left({m\over2}\right)={m\over2}\left({m\over2}+1\right)(1-d)>0$,
we have $M_*(Q_m,d)\ge \lrf{{m\over2}}$.

It remains to show that $M^*(Q_m,d)\le m-1$.
It suffices to prove $d_{m-1}>d_m$,
which is obvious since $d_m=1$ and $d_{m-1}=1+md$.

(ii)\quad
By (i), $m-1$ is the unique mode of $Q_m(x+d)$
if and only if $d_{m-1}>d_{m-2}$.
Note that $d_{m-1}=1+md$ and $d_{m-2}=1+(m-1)d+{m\choose 2}d^2$.
Hence $Q_m(x+d)$ has the unique mode $m-1$ if and only if $0<d<1/{m\choose 2}$.

(iii)\quad
If $0<1-d\le 1/m$, then
$$h\left({m-d\over d+1}\right)={d(m+1)\over(d+1)^2}[3d+1-(1-d)m]>0,$$
which implies that $M_*(Q_m,d)\ge \lrf{{m-d\over d+1}}$.
On the other hand,
$M^*(Q_m,d)\le\overline{m}=\lrc{m-d\over d+1}$ by Corollary \ref{cor2.1}.
Note that
$$\lrf{x}=\left\{
            \begin{array}{ll}
              \lrc{x}, & \hbox{if $x\in{\bf Z}$;} \\
              \lrc{x}-1, & \hbox{otherwise.}
            \end{array}
          \right.$$
Hence $Q_m(x+d)$ has at most two modes $\overline{m}$ and
$\overline{m}-1$, and in particular, only one mode $\overline{m}$ if
${m-d\over d+1}$ is an integer.

(iv)\quad
Denote $t=\lrf{m\over2}$.
Then $M_*(Q_m,d)\ge t$ by (i).
On the other hand, we have by (\ref{dine})
\begin{eqnarray*}
d_{t+1}-d_{t}
&=&\sum\limits_{i=t+1}^{m-1}{i\choose t+1}d^{i-t}-{m\choose t}d^{m-t}\\
&\longrightarrow& \sum\limits_{i=t+1}^{m-1}{i\choose t+1}-{m\choose t}\\
&=&{m\choose t+2}-{m\choose t}
\end{eqnarray*}
when $d$ tends to $1$. Note that ${m\choose t+2}-{m\choose t}<0$.
Hence $d_{t+1}-d_{t}<0$ if $d$ is sufficiently close to $1$, which
implies that $Q_m(x+d)$ has the unique mode $t$. \qed
\begin{rem}\label{rem3.1}
{\em
It is worth pointing out that
modes of $Q_m(x+d)$ may be neither ${\overline m}-1$ nor ${\overline m}$ when $0<d<1$.
For example,
let $1/{m\choose 2}<d<1/m$.
Then ${\overline m}=m$. However,
each mode of $Q_m(x+d)$ is smaller than $m-1$
since $d_{m-2}>d_{m-1}$.
}
\end{rem}

\section{Modes in General Case}
\hspace*{\parindent}
The following theorem shows
the importance of two basic polynomials considered in the last section.
\begin{thm}\label{genpci}
Let $P(x)\in{\bf P}^m_\uparrow$ and $d>0$. Then
$$M_*(Q_m, d)\le M_*(P, d)\le M^*(P, d)\le M^*(x^m, d).$$
Moreover, if $Q_m(x+d)$ has the mode ${\overline m}$,
then so does $P(x+d)$.
In particular,
if $Q_m(x+d)$ has the unique mode ${\overline m}$,
then so does $P(x+d)$ unless $P(x)=x^m$ and
$(m+1)/(d+1)\in{\bf Z}^+$.
\end{thm}
\proof The inequality $M^*(P, d)\le M^*(x^m, d)$ follows from
Corollary \ref{cor2.1} and Proposition \ref{x^m}, so it suffices to
prove the inequality $M_*(Q_m, d)\le M_*(P, d)$.

Let $P(x)=\sum\limits_{j=0}^ma_jx^j$ and
$P(x+d)=\sum\limits_{j=0}^mb_jx^j$.
For $1\le t\le {\overline m}$, let $r=\lrc{(d+1)t}-1$.
Then $t\le r\le m$.
By (\ref{eq1}), we have
\begin{eqnarray*}
td^t(b_{t}-b_{t-1})
&=& \sum\limits_{i=t-1}^ma_id^{i-t}
{i\choose t-1}[(i+1)-(d+1)t]\\
&=&\sum\limits_{i=r}^ma_id^{i}{i\choose t-1}[(i+1)-(d+1)t]\\
&&-\sum\limits_{i=t-1}^{r-1}a_id^{i}
{i\choose t-1}[(d+1)t-(i+1)]\\
&\ge& a_r\sum\limits_{i=r}^md^{i}{i\choose t-1}[(i+1)-(d+1)t]\\ &&-a_r\sum\limits_{i=t-1}^{r-1}
d^{i}{i\choose t-1}[(d+1)t-(i+1)]\\
&=&  a_r\sum\limits_{i=r}^md^{i}{i\choose t-1}[(i+1)-(d+1)t]\\
&=&  a_rtd^t(d_{t}-d_{t-1}),
\end{eqnarray*}
and the equality holds if and only if all $a_i$'s are equal,
i.e., $P$ coincides with $Q_m$.

Take $t=M_*(Q_m,d)$.
Then $d_t>d_{t-1}$ by the definition.
Thus $b_t>b_{t-1}$, which implies that $M_*(P,d)\ge t$,
the desired inequality.

Assume now that ${\overline m}$ is a mode of $Q_m(x+d)$.
Then $d_0\le d_1\le\cdots\le d_{{\overline m}}$.
Thus $b_0\le b_1\le\cdots\le b_{{\overline m}}$.
However,
$b_{{\overline m}}\ge b_{{\overline m}+1}\ge\cdots\ge b_m$ by Corollary \ref{cor2.1}.
Hence ${\overline m}$ is a mode of $P(x+d)$.

In particular, if ${\overline m}$ is the unique mode of $Q_m(x+d)$,
then $M_*(P,d)\ge {\overline m}$. Thus ${\overline m}$ is the unique
mode of $P(x+d)$ if and only if $b_{{\overline m}}>b_{{\overline
m}+1}$, which holds if and only if $P(x)=x^m$ and
$(m+1)/(d+1)\in{\bf Z}^+$ by Corollary \ref{cor2.1} and Proposition
\ref{x^m}. This completes the proof of the theorem. \qed

Combining Theorem \ref{genpci}, Corollary \ref{cor2.1}
and the results of the last section
we conclude that

\begin{cor}
Let $P\in{\bf P}^m_\uparrow$ and $d\ge 1$.
Then $P(x+d)$ has at most two modes ${\overline m}$ and ${\overline m}+1$
if $P(x)=x^m$, or ${\overline m}-1$ and ${\overline m}$ otherwise.
\end{cor}
\begin{cor}\label{P,d=1}
Let $P\in{\bf P}^m_\uparrow$.
Then $P(x+1)$ has the mode $\lrc{m-1\over 2}$.
In particular,
if $P(x)$ is neither $x^m$ nor $\sum\limits_{i=0}^mx^i$,
then $\lrc{m-1\over 2}$ is the unique mode of $P(x+1)$.
\end{cor}
\begin{cor}\label{P,d>1}
Let $d>1$ and $P\in{\bf P}^m_\uparrow$ be such that $P(x)\ne x^m$.
Suppose that one of the following conditions holds:

{\em(i)}\quad
$\overline {m+1}={\overline m}+1$;

{\em(ii)}\quad
${m+1\over d+1}\in{\bf Z}^+$;

{\em(iii)}\quad
$d{\overline m}\in{\bf Z}^+$;

{\em(iv)}\quad
$d\in{\bf Z}^+$;

{\em(v)}\quad
${m\over d+1}\in{\bf Z}^+$.
\\
Then $P(x+d)$ has the unique mode of ${\overline m}$.
\end{cor}

Corollary \ref{P,d=1} and Corollary \ref{P,d>1}(iv) strengthen the
main results of \cite{BM99} and \cite{AABKMR01}, respectively.

In the case $0<d<1$,
the number and location of modes of $P(x+d)$
depend heavily on coefficients of $P(x)$.
Since we are mainly concerned with
those properties of modes
satisfied by generic polynomials in ${\bf P}^m_\uparrow$,
we will not dwell on this case $0<d<1$ any further
but give one useful consequence of
Proposition \ref{Q_m,0<d<1} and Theorem \ref{genpci}, as follows.
\begin{thm}
Let $0<d<1$ and $P\in{\bf P}^m_\uparrow$.
Suppose that $P(x)\ne x^m$.
Then $$\lrf{{m\over2}}\le M_*(P,d)\le M^*(P,d)\le {\overline m}.$$
\end{thm}

\section{Remarks and Open Problems}
\hspace*{\parindent}
Our results can be restated in terms of sequences
instead of polynomials.
For example,
the statement of Conjecture \ref{conj1} is equivalent to the following.
\begin{thm}
Suppose that $0\le a_0\le a_1\le\cdots\le a_ m$ and that $d>0$.
Then the sequence
$$b_j=\sum_{i=j}^ma_id^{i-j}{i\choose j},\quad j=0,1,\ldots,m$$
is unimodal.
\end{thm}

It often occurs that unimodality of a sequence is known,
but to find out the exact number and location of modes of the sequence
is a much more difficult task.
For example,
it is well known that, for each positive integer $n$,
the Stirling number of the second kind $S(n,k)$ is unimodal in $k$
with at most two modes $K_n, K_n+1$,
and that $K_n\sim n/\ln n$.
However it is very difficult to determine
whether the mode of $S(n,k)$ is unique or not.
See \cite{Can78,Har67} for the related results.

We end our paper by proposing the following.
\begin{conj}
Suppose that $P\in{\bf P}^m_\uparrow$
and that $0<d_1<d_2$.
Then $M_*(P,d_1) \ge M_*(P,d_2)$
and $M^*(P,d_1) \ge M^*(P,d_2)$.
\end{conj}

\section*{Acknowledgements}
\hspace*{\parindent} This research was completed during the first
author's stay in the Institute of Mathematics, Academia Sinica,
Taipei. The first author would like to thank the Institute for its
support.

The authors thank the anonymous referees for their valuable
suggestions that led to an improved version of this manuscript.

The first author was partially supported by NSF of Liaoning Province
of China Grant No. 2001102084 and the second author was partially
supported by NSC 92-2115-M-001-016.


\begin{thebibliography}{99}
\bibitem{AABKMR01}
J. Alvarez, M. Amadis, G. Boros, D. Karp, V. H. Moll and L. Rosales,
An extension of a criterion for unimodality,
Electron. J. Combin. 8(2001) \#R30.
\bibitem{BM99}
G. Boros and V. H. Moll, A criterion for unimodality,
Electron. J. Combin. 6(1999) \#R10.
\bibitem{Bre89}
F. Brenti, Unimodal, log-concave, and P\'olya frequency sequences
in combinatorics, Mem. Amer. Math. Soc. 81(1989) no. 413.
\bibitem{Bre94}
F. Brenti,
Log-concave and unimodal sequences in algebra,
combinatorics, and geometry: an update,
Contemp. Math. 178(1994) 71-89.
\bibitem{Can78}E. R. Canfield,
On the location of the maximum Stirling number(s) of the second kind,
Studies in Appl. Math. 59(1978) 83-93.
\bibitem{Har67}L. H. Harper,
Stirling behavior is asymptotically normal,
Ann. Math. Stat. 31(1967) 410-414.
\bibitem{Hog74}
S. G. Hoggar, Chromatic polynomials and logarithmic concavity,
J. Combin. Theory Ser. B 16(1974) 248-254.
\bibitem{Sta89}
R. P. Stanley,
Log-concave and unimodal sequences in algebra,
combinatorics, and geometry,
Ann. New York Acad. Sci. 576(1989) 500-534.
\end{thebibliography}
\end{document}